
\baselineskip=14pt
\parskip=10pt

\magnification=\magstephalf

\def\1{{\overline{1}}}
\def\2{{\overline{2}}}
\parindent=0pt
\overfullrule=0in

\def\frac#1#2{{#1 \over #2}}
\centerline
{\bf Automatic Generation of Convolution Identities for C-finite Sequences}
\bigskip
\centerline
{\it Shalosh B. EKHAD and Doron ZEILBERGER}

\bigskip

{\bf Abstract}: In a recent insightful article, Helmut Prodinger uses sophisticated complex analysis,
with residues, to derive convolution identities for Fibonacci, Tribonacci, and k-bonacci numbers. 
Instead we use a naive-yet-rigorous  `guessing'
approach, using the C-finite ansatz, that can derive such identities in a few seconds, 
but not just for the above-mentioned sequences, but for every C-finite sequence (i.e. a sequence satisfying a linear recurrence with constant coefficients), and even pairs of such sequences.

In a recent delightful article [3], Helmut Prodinger simplified a previous article [2], by Takao Komatsu, that dealt with the binomial convolution
$$
\sum_{k=0}^{n} {{n} \choose {k}} \, T_k T_{n-k} \quad ,
$$

where the $\{T_k\}$  are the {\bf Tribonacci} numbers, that may be defined in terms of their generating function
$$
\sum_{k=0}^{\infty} \, T_k x^k \, = \, \frac{x}{1-x-x^2-x^3} \quad .
$$

Here we suggest an even simpler approach, that is based on {\it guessing}, as preached in [4] (see also [1], Ch. 4).

Recall that a sequence $a(n)$, is $C$-finite  if it satisfies a {\it linear recurrence equation} with {\bf constant} coefficients. Equivalently (see [1][4]) if its
ordinary generating function $\sum_{n=0}^{\infty} a(n)x^n$ is a {\bf rational function} of $x$.
If the order of the recurrence is $d$, then the denominator has degree $d$, and the numerator has degree $\leq d-1$. 
It follows that such a sequence is determined by $2d+1$ constants, and in order to determine them, all 
we need are  $2d+1$ `data points', and use elementary linear algebra, that the computer is glad to do for us.

It is well known and easy to see [1][4] (if the roots of the denominator are distinct) that a $C$-finite sequence of order $d$, let's call it $a(n)$, can ge given by a {\bf Binet}-type formula
$$
a(n) \, = \, \sum_{i=1}^{d} A_i \, \alpha_i^n \quad .
$$

If $b(n)$ is another such sequence, of order $d'$, say, we can write
$$
b(n) \, = \, \sum_{j=1}^{d'} B_j \, \beta_j^n \quad .
$$

Hence the {\bf binomial convolution}
$$
C(n):= \sum_{k=0}^{n} {{n} \choose {k}} \, a(k) b(n-k) \quad 
$$
is a linear combination of powers of the $d\, d'$ numbers
$$
\{ \alpha_i + \beta_j \, | \, 1 \leq i \leq d \, , \, 1 \leq j \leq d' \} \quad,
$$
(why?, exercise left to the reader).  It follows that its generating function is a rational function whose denominator has degree $d\,d'$. Then you  ask the computer to generate the first $2d \,d'+1$ values, and then
to fit that data into the desired rational function (equivalently, recurrence).

If the sequences $a(n)$ and $b(n)$ are the same, then the self-convolution is a linear combination of the powers of the $(d+1)d/2$ numbers
$$
\{ \alpha_i + \alpha_j \, | \, 1 \leq i  \leq j \leq d \} \quad,
$$
and hence the desired rational function only has denominator degree  $(d+1)d/2$, and one only needs $(d+1)d+1$ initial values.

This is implemented in the Maple package {\tt Prodinger.txt} available from 

{\tt https://sites.math.rutgers.edu/\~{}zeilberg/tokhniot/Prodinger.txt} \quad .

At the end of [3], Prodinger stated the generating functions for the  binomial self-convolutions of the tetra-bonacci numbers and the quanta-bonacci numbers, and then added:

{\it ... `but after that the computations become too heavy to be reported here'} \quad .

We filled this {\it much needed gap} and produced all the generating functions of these self-convolutions for $k$-bonacci numbers up to $k=20$. See this output file:

{\tt https://sites.math.rutgers.edu/\~{}zeilberg/tokhniot/oProdinger1.txt} \quad .

We also computed, using our {\it naive-yet-rigorous} approach, the joint  binomial convolution of $k_1$-bonacci and $k_2$-bonacci sequences for $2\leq k_1\leq k_2 \leq 10 $. See the output file

{\tt https://sites.math.rutgers.edu/\~{}zeilberg/tokhniot/oProdinger2.txt} \quad .

More impressively, the Maple package {\tt Prodinger.txt} has the procedures

{\tt ProdingerT(R,x)} and {\tt HelmutT(R1,R2,x)} \quad ,

that can compute, {\bf very fast}, the desired self-convolution or bi-convolution for {\bf any} $C$-finite sequence, or pairs of them, respectively.

We hope that readers will experiment with these two Maple procedures and perhaps find some {\it meta-patterns}. Enjoy!

{\bf References}

[1] Manuel Kauers  and Peter Paule, {\it ``The Concrete Tetrahedron''}, Springer, 2011.

[2] Takao Komatsu. {\it Convolution identities for Tribonacci-type numbers with arbitrary initial values}. Palest. J. Math., {\bf 8(2)} (2019), 413-417.

[3] Helmut Prodinger, {\it Convolutions Identities for Tribonacci numbers via The diagonal of a bivariate generating function},
Palest. J. of Math. {\bf 10}(2) (2021),440-442. Also available here: \hfill\break
{\tt https://arxiv.org/abs/1910.08323} \quad .

[4] Doron Zeilberger, {\it The C-finite Ansatz}, Ramanujan Journal {\bf 31} (2013), 23-32.\hfill\break
{\tt https://sites.math.rutgers.edu/\~{}zeilberg/mamarim/mamarimhtml/cfinite.html} \quad .

\bigskip
\hrule
\bigskip
Shalosh B. Ekhad and Doron Zeilberger, Department of Mathematics, Rutgers University (New Brunswick), Hill Center-Busch Campus, 110 Frelinghuysen
Rd., Piscataway, NJ 08854-8019, USA. \hfill\break
Email: {\tt [ShaloshBEkhad, DoronZeil] at gmail dot com}   \quad .

{\bf Exclusively published in the Personal Journal of Shalosh B. Ekhad and Doron Zeilberger and arxiv.org}

Written: {\bf Aug. 5, 2021}.

\end